# AN IRRATIONALITY MEASURE FOR LIOUVILLE NUMBERS AND CONDITIONAL MEASURES FOR EULER'S CONSTANT

JONATHAN SONDOW

"*Test scores correlate negatively and significantly with the irrationality measure; therefore, students with higher intellectual ability and logical thinking skills tend to be less irrational.*" – Educational research survey

ABSTRACT. The irrationality exponent $\mu(\alpha)$ of an irrational number $\alpha$, defined using the irrationality measure $1/q^\mu$, distinguishes among non-Liouville numbers, and is infinite for Liouville numbers. Using the irrationality measure $1/\beta^q$, we define the *irrationality base* $\beta(\alpha)$, which distinguishes among Liouville numbers and is 1 for non-Liouville numbers. We give some properties and examples. Assuming a condition on certain linear forms in logarithms, for which we present numerical evidence, we prove an upper bound on the irrationality base $\beta(\gamma)$ of Euler's constant, $\gamma$. If $\gamma$ is irrational and the condition turns out to be false in a certain strong sense, we prove an upper bound on the irrationality exponent $\mu(\gamma)$.

## 1. INTRODUCTION

The irrationality exponent $\mu(\alpha)$ of an irrational number $\alpha$ is defined in terms of the irrationality measure $1/q^\mu$ (see Section 2). Using the irrationality measure $1/\beta^q$, we introduce a weaker measure of irrationality, the *irrationality base* $\beta(\alpha)$, as follows. If there exists a real number $\beta \geq 1$ with the property that for any $\varepsilon > 0$ there is a positive integer $q(\varepsilon)$ such that

$$\left|\alpha - \frac{p}{q}\right| > \frac{1}{(\beta+\varepsilon)^q} \quad \text{for all integers } p, q \text{ with } q \geq q(\varepsilon),$$

then we denote by $\beta(\alpha)$ the least such $\beta$, and we call $\beta(\alpha)$ the irrationality base of $\alpha$. If no such $\beta$ exists, then we call $\alpha$ a *super Liouville number,* and we write $\beta(\alpha) = \infty$. Since $\beta(\alpha) = 1$ if $\mu(\alpha)$ is finite (see Lemma 2), *we may regard the irrationality base as a measure of irrationality for Liouville numbers.* We give two examples: a super Liouville number, and a Liouville number with irrationality base 1.

---

2000 *Mathematics Subject Classification*. Primary 11J82, Secondary 11J86.

In [**7**] we gave criteria for irrationality of Euler's constant,

$$\gamma = \lim_{N \to \infty} \left( 1 + \frac{1}{2} + \cdots + \frac{1}{N} - \log N \right).$$

The criteria involve a Beukers-type [**1**] double integral

$$I_n = \iint_{[0,1]^2} -\frac{(x(1-x)y(1-y))^n}{(1-xy)\log xy} \, dx\, dy$$

(defined equivalently in [**8**] as a "hypergeometric integral" or a Nesterenko-type [**5**] series) and a sequence of positive integers

$$S_n = \prod_{m=1}^{n} \prod_{k=0}^{\min(m-1,\, n-m)} \prod_{j=k+1}^{n-k} (n+m)^{\binom{n}{k}^2 \frac{2d_{2n}}{j}},$$

$n = 1, 2, \ldots$, where $d_n$ denotes the least common multiple of the numbers $1, 2, \ldots, n$. In particular, we proved that $\gamma$ is irrational if $(1/n) \log \|\log S_n\|$ does not tend to $-2\log(4/e) = -0.772\ldots$ as $n \to \infty$, where $\|t\|$ denotes the distance from $t$ to the nearest integer. Computations by P. Sebah [**6**] (see Figure 1) suggest that in fact

(S) $$\lim_{n \to \infty} \frac{1}{n} \log \|\log S_n\| = 0.$$

The main result of the present paper is: *If this limit equals zero, then Euler's constant is irrational, but is not a super Liouville number.* More precisely, we prove that *condition* (S) *implies that $\gamma$ has irrationality base $\beta(\gamma) \le 2e = 5.436\ldots$, that is, for any $\varepsilon > 0$ there exists $q(\varepsilon) > 0$ such that*

$$\left| \gamma - \frac{p}{q} \right| > \frac{1}{(2e+\varepsilon)^q} \quad \text{for all integers } p, q, \text{ with } q \ge q(\varepsilon).$$

Thus, a stronger hypothesis than is needed to guarantee irrationality of $\gamma$ yields a stronger conclusion, namely, an irrationality measure for $\gamma$. The proof uses a formula expressing $d_{2n} I_n$ as a $\mathbf{Z}$-linear form in $1, \gamma,$ and $\log S_n$.

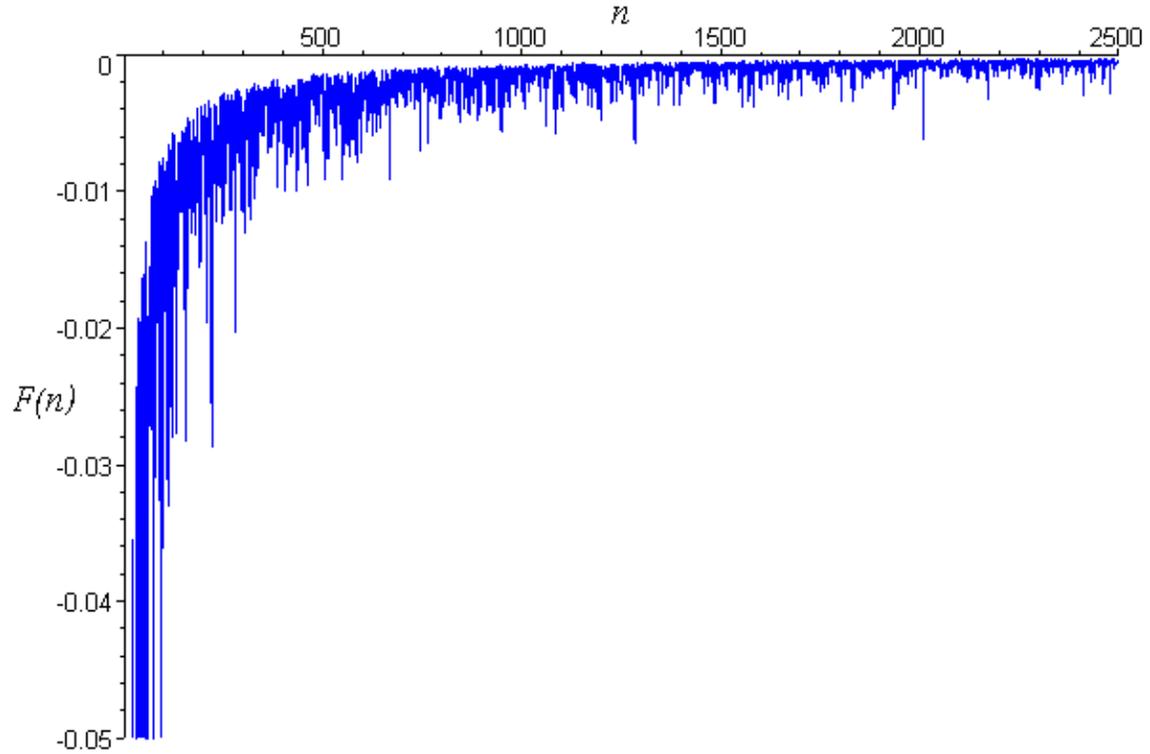

FIGURE 1. Plot of $F(n) = \dfrac{1}{n}\log\|\log S_n\|$

Having given an upper bound for $\beta(\gamma)$ if there does *not* exist a subsequence $\|\log S_{n_k}\|$ tending to zero exponentially, we then give an upper bound for the irrationality exponent $\mu(\gamma)$ if such a subsequence *does* exist, provided both that its convergence rate is different from that of the subsequence $d_{2n_k} I_{n_k}$ (which implies irrationality of $\gamma$), and that the sequence $n_k$ is asymptotically linear. The proof uses standard lemmas on irrationality exponents.

In Section 2, we define irrationality measures, exponents, and bases, and give some properties and examples. In Sections 3 and 4, we prove conditional upper bounds on $\beta(\gamma)$ and $\mu(\gamma)$, respectively. All the conditional bounds are effective.

In a paper in preparation, we construct a sequence of integers similar in size to $S_n$ for which the condition analogous to (S) is false. (This shows that any proof of (S) must use the explicit formula for $S_n$, and cannot rely solely on *a priori* lower bounds for linear forms in logarithms, such as those in Baker's theory.) We also obtain conditional irrationality measures for $\log\pi$, and new ones for $\gamma$ derived from formulas for it in [**9**]. Finally, we give additional properties of the irrationality base, provide examples where it is finite but greater than 1, and fit it into Mahler's classification of transcendental numbers.

## 2. IRRATIONALITY MEASURES, EXPONENTS, AND BASES

**Definition 1.** An *irrationality measure* is a function $f(q,\lambda)$ of a natural number $q$ and a positive real number $\lambda$, which takes values in the positive reals and is decreasing in both $q$ and $\lambda$. Given an irrational number $\alpha$, if there exists $\lambda > 0$ with the property that for any $\varepsilon > 0$ there is a positive integer $q(\varepsilon)$ such that

$$\left|\alpha - \frac{p}{q}\right| > f(q, \lambda+\varepsilon) \quad \text{for all integers } p, q \text{ with } q \geq q(\varepsilon),$$

then we denote by $\lambda(\alpha)$ the least such $\lambda$, and we say that $\alpha$ *has irrationality measure* $f(q, \lambda(\alpha))$.

**Definition 2.** If $\alpha$ has irrationality measure $f(q,\mu) = 1/q^\mu$ for some $\mu = \mu(\alpha)$, then $\mu(\alpha) \in [2, \infty)$ is called the *irrationality exponent* (or, by abuse of terminology, the irrationality measure) of $\alpha$. Otherwise, if no such $\mu$ exists, one writes $\mu(\alpha) = \infty$ and says that $\alpha$ is a *Liouville number*.

**Definition 3.** If $\alpha$ has irrationality measure $f(q,\beta) = 1/\beta^q$, so that $\beta = \beta(\alpha)$ is the least number with the property that for any $\varepsilon > 0$ there exists $q(\varepsilon) > 0$ such that

(1) $$\left|\alpha - \frac{p}{q}\right| > \frac{1}{(\beta+\varepsilon)^q} \quad \text{for all integers } p, q \text{ with } q \geq q(\varepsilon),$$

then we call $\beta(\alpha) \in [1, \infty)$ the *irrationality base* of $\alpha$. Otherwise, if no such $\beta$ exists, we write $\beta(\alpha) = \infty$ and we say that $\alpha$ is a *super Liouville number*.

Note that we may also write inequality (1) as

(2) $$\left|\alpha - \frac{p}{q}\right| > \frac{1}{q^{q \log(\beta+\varepsilon)/\log q}}.$$

**Lemma 1.** (A characterization of super Liouville numbers)
(i). *A real number $\alpha$ has $\beta(\alpha) = \infty$ if and only if, given any $\lambda > 1$, we have*

(3) $$0 < \left|\alpha - \frac{p}{q}\right| < \frac{1}{\lambda^q} \quad \text{for infinitely many integers } p, q, \text{ with } q > 0.$$

(ii). *For fixed $\lambda > 1$, condition* (3) *holds if and only if $\alpha$ is irrational and $\beta(\alpha) \geq \lambda$.*

*Proof.* We prove (ii), which implies (i). Fix $\lambda > 1$. Since $1/\lambda^q < 1/q^\lambda$ for $q$ large, (ii) follows from Definition 2, together with the fact that $\alpha$ is irrational if there exist infinitely many integers $p, q$, with $q > 0$, such that $0 < |\alpha - p/q| < 1/q^\lambda$.   ●

**Example 1.** (*A super Liouville number*) Let $T$ denote the sum of reciprocals of power towers

$$T = \sum_{n=1}^{\infty} \frac{1}{T_n} = \frac{1}{2} + \frac{1}{4^2} + \frac{1}{8^{4^2}} + \cdots,$$

where $T_1 = 2$ and $T_n = (2^n)^{T_{n-1}}$, for $n > 1$. Noting that the denominator of the $n$th partial sum $s_n$ of the series is $T_n$, we see that, given $\lambda > 1$, inequalities (3) hold with $\alpha = T$ and $p/q = s_n$, for all $n$ with $2^n \geq \lambda$. Therefore, $\beta(T) = \infty$.

By the following observation, *the irrationality base measures the irrationality of Liouville numbers*.

**Lemma 2.** (A relation between irrationality exponent and base) *If $\mu(\alpha)$ is finite, then $\beta(\alpha) = 1$; equivalently, if the irrationality base of $\alpha$ exceeds 1, then $\alpha$ is a Liouville number. In particular, a super Liouville number is also a Liouville number.*

*Proof.* This follows from Definitions 2 and 3, since with $\beta = 1$ the exponent $q \log(1+\varepsilon)/\log q$ in (2) tends to infinity with $q$.  •

The converse of Lemma 2 is false, as the following example shows. (For facts used about convergents to continued fractions, see e.g. [**4**, Chapter 10].)

**Example 2.** (*A Liouville number with irrationality base* 1) The continued fraction

$$L = \cfrac{1}{10^{1!} + \cfrac{1}{10^{2!} + \cfrac{1}{10^{3!} + \cdots}}}$$

has irrationality exponent $\mu(L) = \infty$ (see [**4**, page 162, example (b)]). We show that $\beta(L) = 1$. Suppose on the contrary that $\beta(L) > 1$. Then for some $\varepsilon > 0$ there exist infinitely many positive integers $p, q$ such that

(4) $$\left| L - \frac{p}{q} \right| \leq \frac{1}{(1+\varepsilon)^q}.$$

If $q$ is large enough, this bound is less than $1/(2q^2)$, so $p/q$ must be one of the convergents $p_n/q_n$, $n = 0, 1, \ldots$ . Since $q_0 = 1$, $q_1 = 10$, and $q_n = 10^{n!} q_{n-1} + q_{n-2}$ for $n > 1$, by induction it follows that $10^{n!} < q_n < 10^{2n!}$ if $n \geq 2$. Hence, for $n$ sufficiently large,

$$\left| L - \frac{p_n}{q_n} \right| > \frac{1}{q_n q_{n+2}} > \frac{1}{10^{2(n!+(n+2)!)}} > \frac{1}{(1+\varepsilon)^{10^{n!}}} > \frac{1}{(1+\varepsilon)^{q_n}},$$

contradicting (4). Therefore, $\beta(L) = 1$.

In Section 3, assuming a hypothesis, we deduce upper bounds for an irrationality base directly from the definition. In Section 4, however, in order to obtain conditional upper bounds for the irrationality exponent of a real number $\alpha$, we first construct sequences of integers $p_k, q_k$ such that

$$\varepsilon_k := q_k \alpha - p_k$$

tends to zero exponentially. We then use one of the following two standard lemmas. (For proofs, see [**2**, Lemma 3.5] and [**3**, Remark 2.1], respectively.)

**Lemma 3.** (Chudnovsky) *Suppose that*

$$\limsup_{k \to \infty} \frac{1}{k} \log |q_k| \leq \sigma, \quad \lim_{k \to \infty} \frac{1}{k} \log |\varepsilon_k| = -\tau,$$

*for some positive numbers $\sigma, \tau$. Then $\alpha$ has irrationality exponent $\mu(\alpha) \leq 1 + \sigma/\tau$.*

**Lemma 4.** (Hata) *Suppose that $\alpha$ is irrational and that*

$$\lim_{k \to \infty} \frac{1}{k} \log |q_k| = \sigma, \quad \limsup_{k \to \infty} \frac{1}{k} \log |\varepsilon_k| \leq -\tau,$$

*for some positive numbers $\sigma, \tau$. Then $\mu(\alpha) \leq 1 + \sigma/\tau$.*

**Remark 1.** If $|q_k| \sim e^{\sigma k}$ and $|\varepsilon_k| \sim e^{-\tau k}$ as $k \to \infty$, then $|\alpha - p_k/q_k| \sim |q_k|^{-1-\tau/\sigma}$. By Definition 2, it follows that $\mu(\alpha) \geq 1 + \tau/\sigma$. Hence $\sigma \geq \tau$, by Lemma 3.

**Remark 2.** It would be interesting to find analogous lemmas for the irrationality base.

## 3. CONDITIONAL BOUNDS ON THE IRRATIONALITY BASE OF EULER'S CONSTANT

Assuming Condition (S) or weaker conditions, we derive upper bounds on the irrationality base of $\gamma$. Note that (S) implies that $(1/n) \log \{\log S_n\}$ also tends to zero, where $\{t\}$ denotes the fractional part of $t > 0$.

From [**7**], [**8**], we have the relations

(5) $$d_{2n} I_n = d_{2n} \binom{2n}{n} \gamma + \log S_n - d_{2n} A_n, \quad d_{2n} A_n \in \mathbf{Z}$$

7and the criterion

(6) $$d_{2n} I_n = \{\log S_n\}, \text{ for } n \geq n_0 \Leftrightarrow \gamma \in \mathbf{Q}.$$

From [**7**], Stirling's formula, and the Prime Number Theorem, we have

(7) $$4^{-2n} > I_n = 4^{-2n(1+o(1))},$$

(8) $$2^{2n} > \binom{2n}{n} = 2^{2n(1+o(1))},$$

(9) $$e^{n(1+\varepsilon)} > d_n = e^{n(1+o(1))},$$

respectively, where the first inequality in (9) holds for any $\varepsilon > 0$ and $n \geq n(\varepsilon)$.

Criterion (6) and asymptotics (7), (9) yield the implication

$$\gamma \in \mathbf{Q} \Rightarrow \lim_{n \to \infty} \frac{1}{n} \log\{\log S_n\} = -2\log(4/e).$$

Each of the conditions that we assume below contradicts the last equation, and therefore implies that $\gamma \notin \mathbf{Q}$.

To find a conditional upper bound for $\beta(\gamma)$, suppose first that $\beta(\gamma) > 1$. It follows from Definition 3 that, given $\lambda$ between 1 and $\beta(\gamma)$, there exist integers $p_k, q_k$, for $k = 1, 2, \ldots$, such that $0 < q_k < q_{k+1}$ and

(10) $$\left| \gamma - \frac{p_k}{q_k} \right| < \frac{1}{\lambda^{q_k}}.$$

Set

$$n_k = \left\lceil \frac{q_k}{2} \right\rceil$$

and use (5) to write

(11) $$d_{2n_k}\left(I_{n_k} - \binom{2n_k}{n_k}\left(\gamma - \frac{p_k}{q_k}\right)\right) = \log S_{n_k} - d_{2n_k}\left(A_{n_k} - \binom{2n_k}{n_k}\frac{p_k}{q_k}\right),$$

where $d_{2n_k} A_{n_k} \in \mathbf{Z}$. Since $q_k \leq 2n_k$, the expression following $\log S_{n_k}$ is an integer; thus the absolute value of the left-hand side of (11), if less than 1, is $\geq \|\log S_{n_k}\|$. Also, since $q_k \geq 2n_k - 1$, inequalities (7), (8), (9), (10) imply that, for any $\varepsilon > 0$, we have

$$(12) \quad d_{2n_k} \left| I_{n_k} - \binom{2n_k}{n_k}\left(\gamma - \frac{p_k}{q_k}\right) \right| < e^{2n_k(1+\varepsilon)}\left(4^{-2n_k} + \lambda(\lambda/2)^{-2n_k}\right) \quad \text{for } k \geq k(\varepsilon).$$

Suppose now that $\beta(\gamma) > 2e$. Choose $\lambda$ between $2e$ and $\beta(\gamma)$, and let $\varepsilon > 0$ be so small that

$$(13) \quad m_{\lambda,\varepsilon} := e^{-1-\varepsilon}\min(4, \lambda/2) > 1.$$

Then, for $k$ large, the right-hand side of (12) is less than 1. Using (11), (12), it follows that

$$(14) \quad \left\|\log S_{n_k}\right\| < (1+\lambda)m_{\lambda,\varepsilon}^{-2n_k} \quad \text{for } k \geq k(\lambda,\varepsilon),$$

which contradicts (S). This establishes our main result.

**Theorem 1.** (First conditional bound on $\beta(\gamma)$) *Assume that*

$$\lim_{n\to\infty}\frac{1}{n}\log\left\|\log S_n\right\| = 0.$$

*Then Euler's constant is irrational, but is not a super Liouville number. More precisely, the irrationality base of $\gamma$ satisfies $\beta(\gamma) \leq 2e$, that is, for any $\varepsilon > 0$ there exists $q(\varepsilon) > 0$ such that*

$$\left|\gamma - \frac{p}{q}\right| > \frac{1}{(2e+\varepsilon)^q} \quad \text{for all integers } p, q \text{ with } q \geq q(\varepsilon).$$

Letting $(\lambda,\varepsilon)$ tend to $(\beta(\gamma),0)$ in (14), and using (13), we obtain

$$\liminf_{n\to\infty}\frac{1}{n}\log\left\|\log S_n\right\| \leq -2\log\min(4/e, \beta(\gamma)/(2e)).$$

This yields the following generalization of Theorem 1, giving an upper bound on $\beta(\gamma)$ assuming a condition possibly less difficult to verify than (S).

**Theorem 2.** (Second conditional bound on $\beta(\gamma)$) *If*

$$(15) \quad \liminf_{n\to\infty}\frac{1}{n}\log\left\|\log S_n\right\| \geq -\delta$$

*for some non-negative number $\delta < 2\log(4/e)$, then $\beta(\gamma) \leq 2e^{1+(\delta/2)} < 8$.*

**Remark.** For $1 \leq n \leq 2500$, the minimum of $(1/n) \log \|\log S_n\|$ is $-0.667\ldots$ (at $n = 5$), except for the value $-1.480\ldots$ at $n = 1$. Figure 1 suggests that (15) holds with $\delta = 0.01$, and probably also with $\delta = 0$, which is Condition (S).

Finally, suppose that $8 < \lambda < \beta(\gamma)$. We deduce from (10), (14) and asymptotics (8), (9) that, for $k$ large, the left-hand side of (11) is positive and less than 1, hence equal to $\{\log S_{n_k}\}$. Since now $m_{\lambda,\varepsilon} = 4e^{-1-\varepsilon}$, we obtain

$$\{\log S_{n_k}\} < (1+\lambda)\,(4/e^{1+\varepsilon})^{-2n_k},$$

for $k \geq k(\lambda, \varepsilon)$. Letting $\varepsilon$ tend to zero, we obtain the following result, which assumes a weaker condition than the previous ones.

**Theorem 3.** (Third conditional bound on $\beta(\gamma)$) *We have* $\beta(\gamma) \leq 8$ *if*

$$\liminf_{n \to \infty} \frac{1}{n} \log\{\log S_n\} > -2\log(4/e).$$

## 4. CONDITIONAL BOUNDS ON THE IRRATIONALITY EXPONENT OF EULER'S CONSTANT

We have shown that $\gamma$ has irrationality measure $1/\beta^q$, with $\beta = \beta(\gamma) < \infty$, if there do *not* exist subsequences $\|\log S_{n_k}\|$ tending to zero exponentially. We now show that $\gamma$ has irrationality measure $1/q^\mu$, with $\mu = \mu(\gamma) < \infty$, if such a subsequence *does* exist, provided both that its convergence rate is different from that of $d_{2n_k} I_{n_k}$, and that the sequence $n_k$ is asymptotically linear. More generally, the limit of the subsequence may be any rational number $a/b$ in the interval $[0,1]$, and $\|\log S_{n_k}\|$ may be replaced by $\{\log S_{n_k}\}$.

**Theorem 4.** (Conditional bounds on $\mu(\gamma)$) *Assume that there exists a sequence of positive integers $n_k$, $k = 1, 2, \ldots$, such that*

$$(16)_\sigma,\ (16)_\tau \qquad \lim_{k \to \infty} \frac{n_k}{k} = \sigma, \qquad \lim_{k \to \infty} \frac{1}{k} \log\left|\{\log S_{n_k}\} - \frac{a}{b}\right| = -\tau,$$

*for some integers $a, b$ and positive numbers $\sigma, \tau$, with $\tau \neq 2\sigma \log(4/e)$. Then $\gamma$ has irrationality exponent $\mu(\gamma) \leq \mu_{\sigma,\tau}$, where*

$$(17) \quad \mu_{\sigma,\tau} = \begin{cases} 1 + \dfrac{2\sigma}{\tau}\log(2e) & \text{if } \tau < 2\sigma\log(4/e), \\ \dfrac{\log 8}{\log 4 - 1} = 5.383\ldots & \text{if } \tau > 2\sigma\log(4/e), \end{cases}$$

*that is, for any $\varepsilon > 0$ there exists $q(\varepsilon) > 0$ such that*

$$\left|\gamma - \frac{p}{q}\right| > \frac{1}{q^{\mu_{\sigma,\tau}+\varepsilon}} \quad \text{for all integers } p, q, \text{ with } q \geq q(\varepsilon).$$

*In particular, the hypothesis implies that $\gamma$ is irrational, but is not a Liouville number. The same result holds with $\{\log S_{n_k}\}$ replaced by $\|\log S_{n_k}\|$.*

*Proof.* Define integers $p_k, q_k$, for $k = 1, 2, \ldots$, by the formulas

$$(18) \quad p_k = b\, d_{2n_k} A_{n_k} - b\lfloor \log S_{n_k} \rfloor - a, \qquad q_k = b\, d_{2n_k} \binom{2n_k}{n_k}.$$

According to (5), we have

$$(19) \quad b^{-1}(q_k \gamma - p_k) = d_{2n_k} I_{n_k} + \frac{a}{b} - \{\log S_{n_k}\}.$$

Asymptotics (7), (8), (9), together with limit $(16)_\sigma$, imply that

$$\lim_{k\to\infty} \frac{1}{k} \log |q_k| = 2\sigma \log(2e)$$

and

$$(20) \quad \lim_{k\to\infty} \frac{1}{k} \log(d_{2n_k} I_{n_k}) = -2\sigma \log(4/e).$$

Using $\tau \neq 2\sigma \log(4/e)$, we deduce from (19), (20), $(16)_\tau$ that

$$\lim_{k\to\infty} \frac{1}{k} \log |q_k \gamma - p_k| = -\min(\tau, 2\sigma \log(4/e)).$$

Lemma 3 now implies the theorem, except for its last assertion, whose proof requires the following modifications. In $(16)_\tau$, replace $\{\log S_{n_k}\}$ by $\|\log S_{n_k}\|$. If the two numbers are equal, define $p_k$ as in (18); if they are unequal, define $p_k$ by

$$p_k = b\, d_{2n_k} A_{n_k} - b\lceil \log S_{n_k} \rceil + a.$$

In both cases, define $q_k$ as in (18). Then (19) becomes

$$b^{-1}(q_k \gamma - p_k) = d_{2n_k} I_{n_k} \pm \left( \| \log S_{n_k} \| - \frac{a}{b} \right)$$

and the rest of the proof goes through unchanged. ●

Using Lemma 4 in place of Lemma 3, we can deduce the conclusions of Theorem 4.1 from a weaker hypothesis, as follows.

**Theorem 5.** *Assume that there exists a sequence $n_k$ such that* $(16)_\sigma$ *and*

(21) $$\limsup_{k \to \infty} \frac{1}{k} \log \left| \{\log S_{n_k}\} - \frac{a}{b} \right| \leq -\tau$$

*hold, for some integers $a, b$ and positive numbers $\sigma, \tau$. Assume further that $-2\sigma \log(4/e)$ is not the limit of any subsequence of $(1/k) \log |\{\log S_{n_k}\} - a/b|$. Then $\gamma$ has irrationality exponent $\mu(\gamma) \leq \mu_{\sigma,\tau}$, where $\mu_{\sigma,\tau}$ is given by (17). The same result holds with $\{\log S_{n_k}\}$ replaced by $\|\log S_{n_k}\|$.*

*Proof.* We repeat the previous proof through (20). Using the subsequence condition, we deduce from (19), (20), (21) that

$$\limsup_{k \to \infty} \frac{1}{k} \log |q_k \gamma - p_k| \leq -\min(\tau, 2\sigma \log(4/e)).$$

By Lemma 4, it only remains to show that, under the hypothesis, $\gamma$ is irrational. Suppose on the contrary that $\gamma \in \mathbf{Q}$. Then, by criterion (6) and inequalities (7), (9), we have

(22) $$\{\log S_{n_k}\} = d_{2n_k} I_{n_k} \to 0 \quad \text{as } k \to \infty.$$

Hence $a/b = 0$, by (21). But then (20) and the equality in (22) contradict the subsequence hypothesis. Thus, $\gamma \notin \mathbf{Q}$ and the proof for $\{\log S_{n_k}\}$ is complete. The same modification as above works for $\|\log S_{n_k}\|$. ●

**Remark.** Alternatively, one could prove Theorem 5 first and obtain Theorem 4 as a corollary. We chose the present order because Theorem 4 has a simpler statement and proof than Theorem 5.

**ACKNOWLEDGEMENTS.** The author is grateful to T. Rivoal for asking about conditional measures for $\gamma$ and for suggesting part of Theorem 4, to P. Sebah for Figure 1 and computations of $\|\log S_n\|$, and to W. Zudilin for advice on the definitions and terminology in Section 2.

## REFERENCES

1. F. Beukers, A note on the irrationality of $\zeta(2)$ and $\zeta(3)$, *Bull. London Math. Soc.* **12** (1979) 268-272.
2. G. Chudnovsky, Hermite-Padé approximations to exponential functions and elementary estimates of the measure of irrationality of $\pi$, *Lecture Notes in Math.* **925** (1982) 299-325.
3. M. Hata, Rational approximations to $\pi$ and some other numbers, *Acta Arith.* **63** (1993) 335-349.
4. G. H. Hardy and E. M. Wright, *An Introduction to the Theory of Numbers*, 5th edition, Oxford Univ. Press, New York, 1979.
5. Y. Nesterenko, A few remarks on $\zeta(3)$, *Math. Notes* **59** (1996) 625-636.
6. P. Sebah, personal communication, 29 August 2002.
7. J. Sondow, Criteria for irrationality of Euler's constant, *Proc. Amer. Math. Soc.,* to appear.
8. \_\_\_\_\_\_\_\_, A hypergeometric approach, via linear forms involving logarithms, to irrationality criteria for Euler's constant, *CRM Conf. Proc. of CNTA* 7 (*May, 2002*), submitted.
9. \_\_\_\_\_\_\_\_ and W. Zudilin, Euler's constant, $q$-logarithms, and formulas of Ramanujan and Gosper, *Ramanujan J.*, submitted.

209 WEST 97TH STREET, NEW YORK, NY 11025
*email address*: jsondow@alumni.princeton.edu